\documentclass[11pt,reqno]{amsart}

\usepackage{amsmath,amssymb,amsfonts,graphicx,hyperref,subcaption}
\usepackage{mathtools}
\usepackage{geometry}
\title[Inverse Chafee--Infante via WGAN-GP]{Backward Reconstruction of the Chafee--Infante Equation via Physics-Informed WGAN-GP}

\author{Joseph L. Shomberg}
\address{Department of Mathematics and Computer Science\\
Providence College\\
Providence, RI 02918, USA}
\email{jshomber@providence.edu}

\date{\today}

\begin{document}

\maketitle

\begin{abstract}
We present a physics-informed Wasserstein GAN with gradient penalty (WGAN-GP) for solving the inverse Chafee--Infante problem on two-dimensional domains with Dirichlet boundary conditions. 
The objective is to reconstruct an unknown initial condition from a near-equilibrium state obtained after 100 explicit forward Euler iterations of the reaction-diffusion equation
\[
u_t - \gamma\Delta u + \kappa\left(u^3 - u\right)=0.
\]
Because this mapping strongly damps high-frequency content, the inverse problem is severely ill-posed and sensitive to noise.

Our approach integrates a U-Net generator, a PatchGAN critic with spectral normalization, Wasserstein loss with gradient penalty, and several physics-informed auxiliary terms, including Lyapunov energy matching, distributional statistics, and a crucial forward-simulation penalty. 
This penalty enforces consistency between the predicted initial condition and its forward evolution under the \emph{same} forward Euler discretization used for dataset generation. 
Earlier experiments employing an Eyre-type semi-implicit solver were not compatible with this residual mechanism due to the cost and instability of Newton iterations within batched GPU training.

On a dataset of 50k training and 10k testing pairs on $128\times128$ grids (with natural $[-1,1]$ amplitude scaling), the best trained model attains a mean absolute error (MAE) of approximately \textbf{0.23988159} on the full test set, with a sample-wise standard deviation of about \textbf{0.00266345}. 
The results demonstrate stable inversion, accurate recovery of interfacial structure, and robustness to high-frequency noise in the initial data.
\end{abstract}

%%%%%%%%%%%%%%%%%%%%%%%%%%%%%%%%%%%%%%%%
\section{Introduction}
%%%%%%%%%%%%%%%%%%%%%%%%%%%%%%%%%%%%%%%%

Inverse problems for nonlinear parabolic PDEs are severely ill-posed: small perturbations in observed data may correspond to large deviations in the underlying initial condition. 
The Chafee--Infante equation is a prototypical bistable reaction-diffusion model whose forward dynamics rapidly relax toward the wells $\pm1$, erasing high-frequency information and producing sharp interfacial patterns. 
Reconstructing the initial state from a deeply relaxed near-equilibrium snapshot is therefore a canonical and challenging inverse problem.

Classical approaches to initial data reconstruction for nonlinear parabolic equations appears in \cite{Tikhonov,IsakovInverse}, but standard regularization techniques tend to oversmooth the recovered solution and fail to reproduce sharp interface geometry \cite{RudinOsherFatemi}. 
PDE-constrained optimization methods \cite{PDEControl} are computationally expensive and remain sensitive to the underlying ill-posedness.

Recent advances in generative modeling \cite{GAN,WGAN,WGAN-GP} and physics-informed machine learning \cite{PINN,LiPDEML} offer an alternative pathway. 
By combining adversarial learning with structural constraints derived from the governing PDE, physics-informed GANs can represent multimodal, nonconvex solution distributions while enforcing physical consistency. 
In this work, we demonstrate that a WGAN-GP augmented with a forward-simulation residual term can stably reconstruct high-frequency initial conditions of the Chafee--Infante equation from near-equilibrium data.

%%%%%%%%%%%%%%%%%%%%%%%%%%%%%%%%%%%%%%%%
\section{Forward Dataset Generation}
%%%%%%%%%%%%%%%%%%%%%%%%%%%%%%%%%%%%%%%%

Training and testing data are constructed using an explicit forward Euler 
scheme for the Chafee--Infante equation
\[
u_t = \gamma\Delta u - \kappa\left(u^3-u\right), \qquad \gamma=0.005,\quad \kappa=4.7,
\]
on the domain $\Omega=[-1,1]^2$ with Dirichlet boundary conditions, discretized 
on a $128\times128$ Cartesian grid. 

Random initial conditions are drawn as
\[
u_0(x) \sim \mathrm{Unif}([-0.02,0.02]) \quad \text{for } x \in \Omega,
\]
with Dirichlet boundary conditions imposed on $\partial\Omega$. 
This choice introduces small-scale, high-frequency perturbations in the interior
while remaining naturally compatible with the bistable wells $\pm 1$, so no
additional amplitude normalization is required.
Each initial condition is evolved for 100 forward Euler time steps with 
timestep $\Delta t=10^{-3}$, producing a near-equilibrium field $u_{100}$. 
We store the resulting datasets as \emph{source-target} pairs
\[
\left(\text{src},\text{tar}\right) = \left(u_{100},u_0\right).
\]

Earlier stages of this project employed an Eyre-type semi-implicit convex-splitting scheme with Newton iterations. 
While unconditionally stable for forward simulation, such solvers are incompatible with the residual-based training loop used here: each time step requires an inner nonlinear solver whose control flow and convergence behavior are ill-suited to batched, differentiable GPU execution. 
For this reason, the present study uses forward Euler exclusively for both dataset generation and residual evaluation, ensuring numerical consistency.

%%%%%%%%%%%%%%%%%%%%%%%%%%%%%%%%%%%%%%%%
\section{Model Architecture and Physics-Informed Losses}
%%%%%%%%%%%%%%%%%%%%%%%%%%%%%%%%%%%%%%%%

%%%%%%%%%%%%%%%%%%%%%%%%%%%%%%%%%%%%%%%%
\subsection{Generator $+$ Critic}
%%%%%%%%%%%%%%%%%%%%%%%%%%%%%%%%%%%%%%%%
The generator is a U-Net with skip connections and a final $\tanh$ activation 
to enforce values in $[-1,1]$. 
Instance normalization is used throughout the network. 
This architecture is well-suited for multiscale pattern inversion.

A PatchGAN-style critic with spectral normalization evaluates local consistency 
between $(u_{100}, u_0)$ and $(u_{100}, \hat{u}_0)$, where $u_{100}$ is a given source image, $u_{0}$ is the ground truth/target initial condition, and $\hat{u}_0$ is generated/fake. 
We employ the Wasserstein loss with gradient penalty.

Our generator loss is as follows:
\[
\mathcal{L}_G = 
- \mathbb{E}[D(u_{100},\hat{u}_0)] 
+ \lambda_E \mathcal{L}_{\text{energy}}
+ \lambda_R \mathcal{L}_{\text{res}}
+ \lambda_{\text{MAE}}\!\cdot\!\text{MAE}
+ \lambda_{\mu}\mathcal{L}_{\text{mean}}
+ \lambda_{\sigma}\mathcal{L}_{\text{var}}.
\]

\paragraph{\textbf{Lyapunov Energy Loss.}}
We compute the discrete free energy
\[
\mathcal{E}(u)=\int \left(\frac{\gamma}{2}|\nabla u|^2+\frac{\kappa}{4}\left(u^2-1\right)^2\right)dx,
\]
and penalize mismatches between $\mathcal{E}(\hat{u}_0)$ and $\mathcal{E}(u_0)$.

\paragraph{\textbf{Statistical Moment Matching.}}
We match means and variances between predicted and true initial data.  
This stabilizes global amplitude and symmetry.

\paragraph{\textbf{Residual (Forward-Simulation Consistency).}}
This is the central contribution.  
After predicting $\hat{u}_0$, we run a forward Euler simulation for $s$ steps:
\[
F^s\left(\hat{u}_0\right) \approx u_{\text{sim}}.
\]
We then compare $u_{\text{sim}}$ with the observed $u_{100}$,
\[
\mathcal{L}_{\text{res}}
= \left\|F^s\left(\hat{u}_0\right) - u_{100}\right\|_1.
\]

Using the same forward Euler discretization for both data generation and residual evaluation ensures numerical consistency and avoids the instability and computational cost associated with Newton-based solvers.
Using forward Euler exclusively avoids these issues and yields a clean, self-consistent 
residual term.

%%%%%%%%%%%%%%%%%%%%%%%%%%%%%%%%%%%%%%%%
\section{Training and Optimization}
%%%%%%%%%%%%%%%%%%%%%%%%%%%%%%%%%%%%%%%%

Training used batch size $1$, Adam optimizers with $\beta_1=0.5$, and gradient penalty coefficient $\lambda_{GP}=10$. 
Models were trained for several thousand iterations with checkpointing and automatic selection based on validation MAE. 
The weights $(\lambda_E,\lambda_R,\lambda_{\text{MAE}})$ were chosen 
heuristically via informed trial and error, with limited Bayesian optimization 
used only as guidance.
 
%%%%%%%%%%%%%%%%%%%%%%%%%%%%%%%%%%%%%%%%
\section{Results \& Figures}
%%%%%%%%%%%%%%%%%%%%%%%%%%%%%%%%%%%%%%%%

During training, model selection was guided by a fixed three-sample validation MAE. 
After training, performance was evaluated on the full test set of 10{,}000 samples.

On 10{,}000 unseen test pairs, the best model attains
\[
\text{MAE}_{\text{test}} \approx 0.23988159 \pm 0.00266345,
\]
where the $\pm$ denotes the standard deviation of the sample-wise MAE over the test set. 
The corresponding standard error of the mean is
\[
\mathrm{SEM} \approx \frac{0.00266345}{\sqrt{10{,}000}} \approx 2.663\times 10^{-5},
\]
so the reported mean MAE is statistically well resolved.

Qualitatively, the generator recovers the correct mixture of broad domains 
and thin transition layers characteristic of bistable dynamics. 
The residual term ensures that predicted initial conditions evolve 
accurately into the observed near-equilibrium patterns under the 
forward Euler dynamics.

\begin{figure}[h]
    \centering
    % Row 1
    \begin{subfigure}[b]{0.3\textwidth}
        \includegraphics[width=\textwidth]{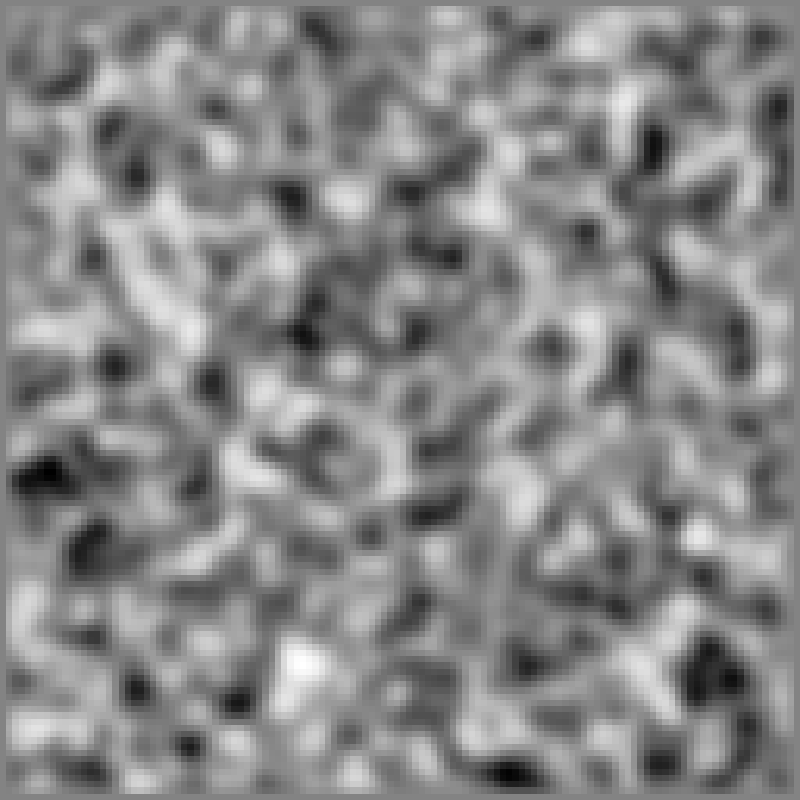}
        \caption*{Source image 1}
    \end{subfigure}
    \hfill
    \begin{subfigure}[b]{0.3\textwidth}
        \includegraphics[width=\textwidth]{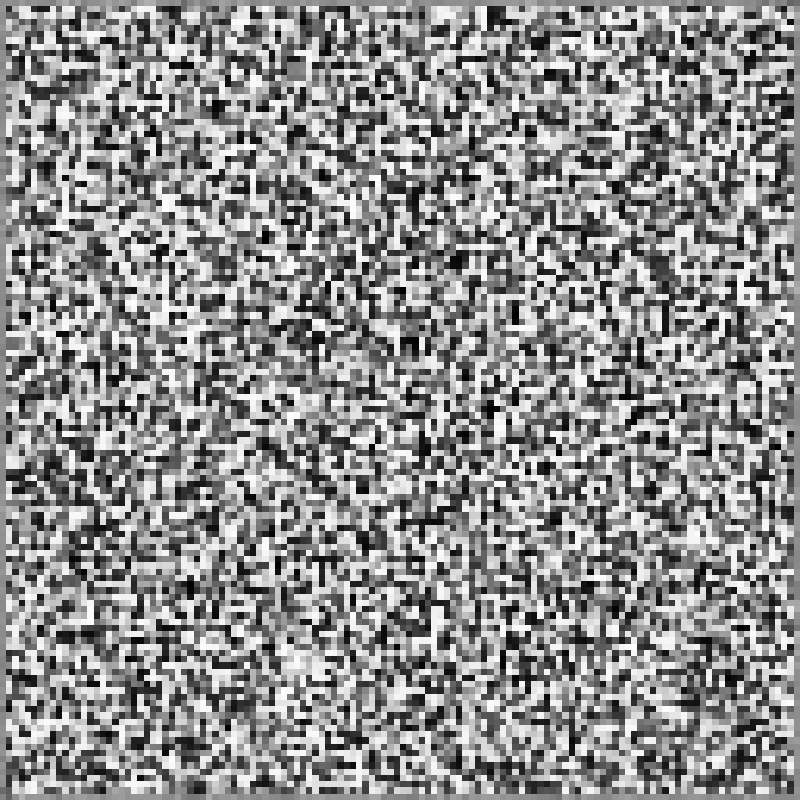}
        \caption*{Generated image 1}
    \end{subfigure}
    \hfill
    \begin{subfigure}[b]{0.3\textwidth}
        \includegraphics[width=\textwidth]{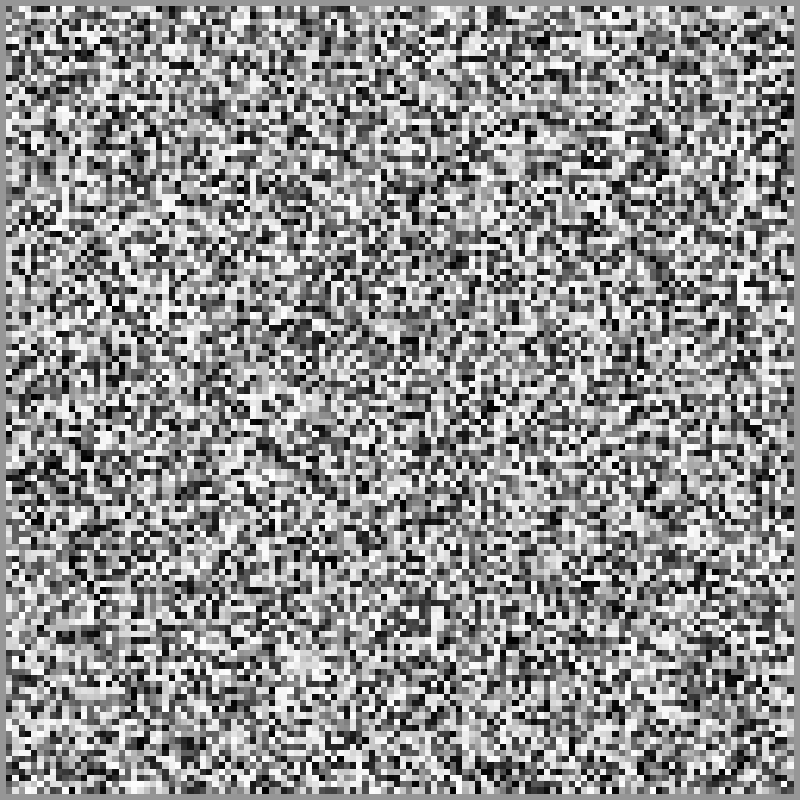}
        \caption*{Target image 1}
    \end{subfigure}

    \vspace{1em} % Adds vertical space between rows

    % Row 2
    \begin{subfigure}[b]{0.3\textwidth}
        \includegraphics[width=\textwidth]{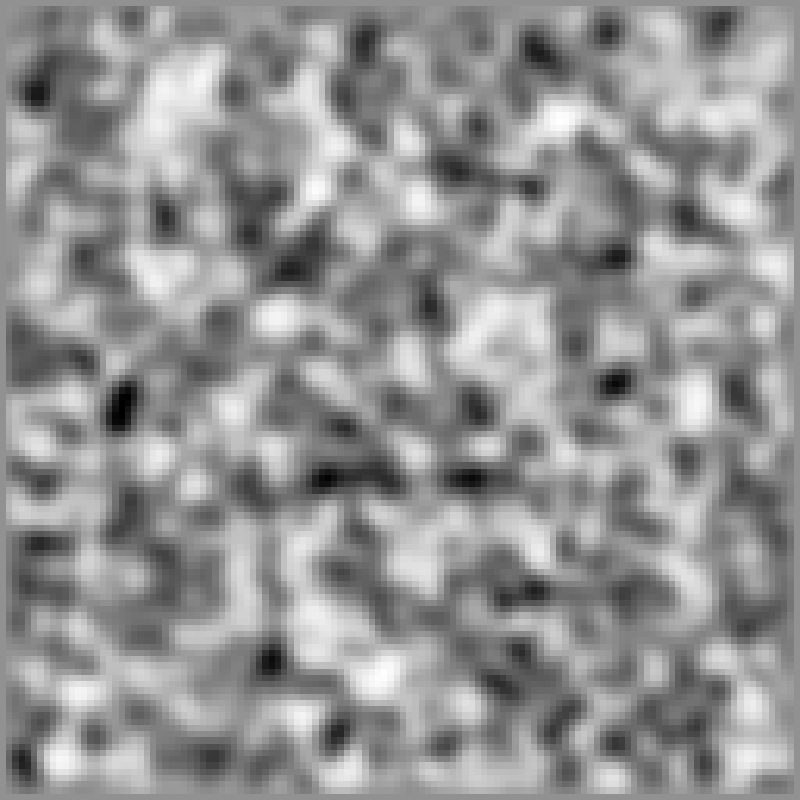}
        \caption*{Source image 2}
    \end{subfigure}
    \hfill
    \begin{subfigure}[b]{0.3\textwidth}
        \includegraphics[width=\textwidth]{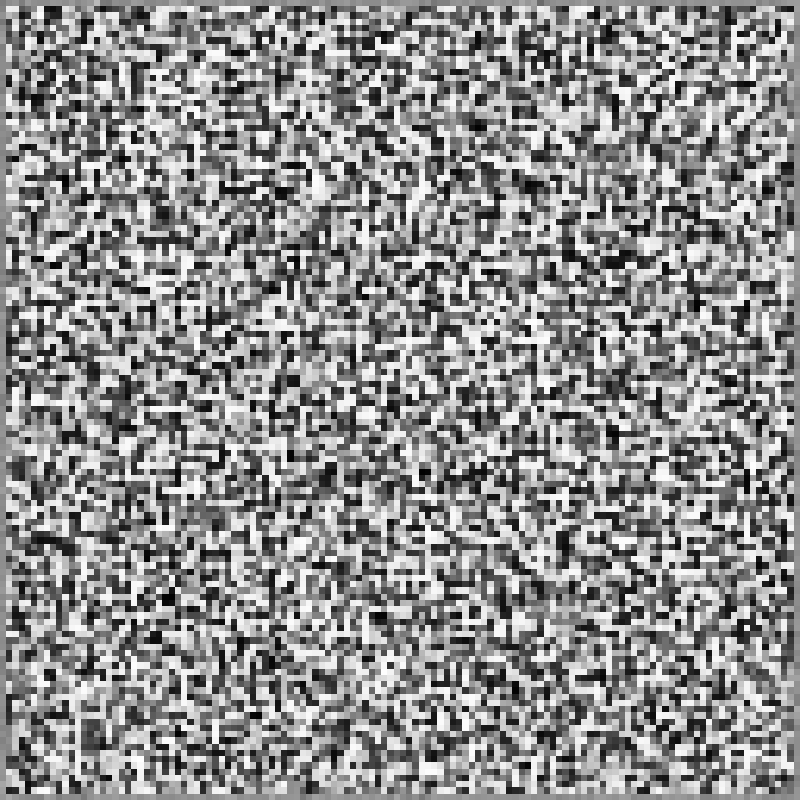}
        \caption*{Generated image 2}
    \end{subfigure}
    \hfill
    \begin{subfigure}[b]{0.3\textwidth}
        \includegraphics[width=\textwidth]{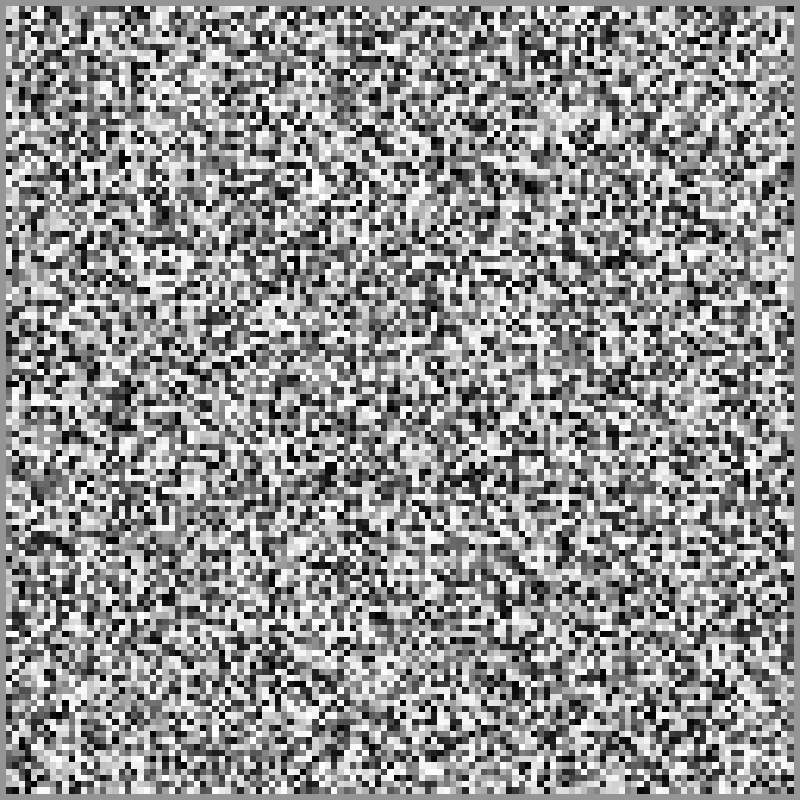}
        \caption*{Target image 2}
    \end{subfigure}

    \vspace{1em} % Adds vertical space between rows

    % Row 3
    \begin{subfigure}[b]{0.3\textwidth}
        \includegraphics[width=\textwidth]{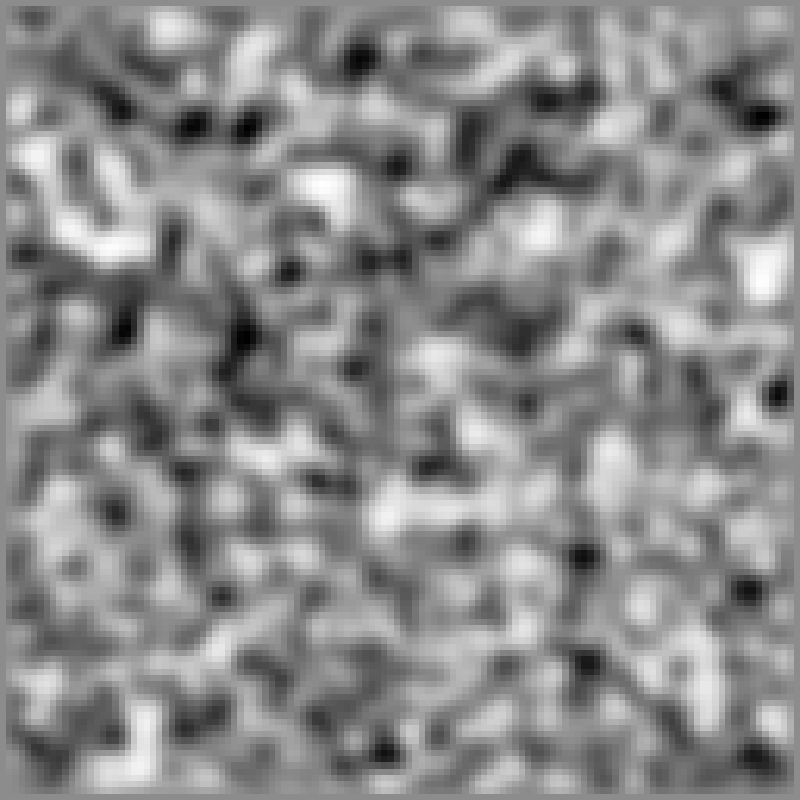}
        \caption*{Source image 3}
    \end{subfigure}
    \hfill
    \begin{subfigure}[b]{0.3\textwidth}
        \includegraphics[width=\textwidth]{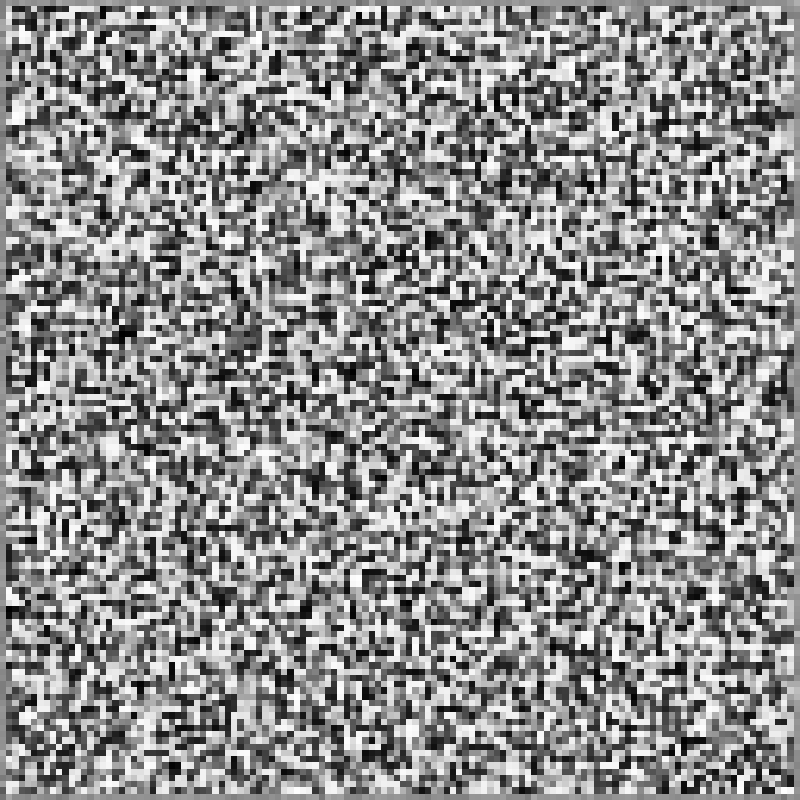}
        \caption*{Generated image 3}
    \end{subfigure}
    \hfill
    \begin{subfigure}[b]{0.3\textwidth}
        \includegraphics[width=\textwidth]{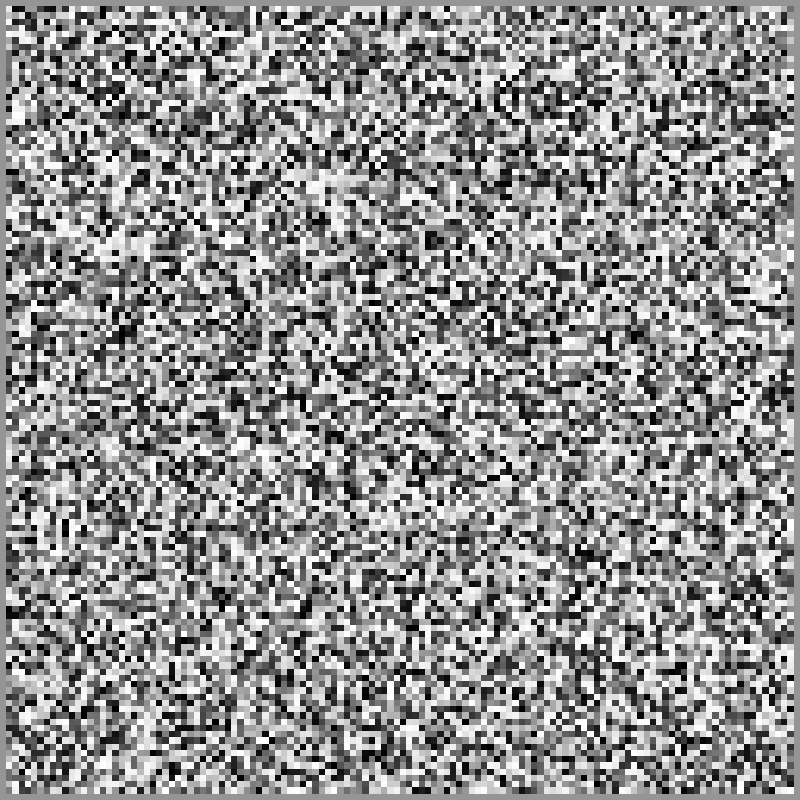}
        \caption*{Target image 3}
    \end{subfigure}

    \caption{MAE Testing: $model(source)=generated$ which is compared to $target$.}
\end{figure}

%%%%%%%%%%%%%%%%%%%%%%%%%%%%%%%%%%%%%%%%
\section{Discussion}
%%%%%%%%%%%%%%%%%%%%%%%%%%%%%%%%%%%%%%%%

The central finding is that physics-informed generative inversion of dissipative PDEs succeeds only when the physics and numerics used inside the loss match those that generated the data. By committing to forward Euler for both dataset generation and residual evaluation, the residual enforces genuine discrete dynamical consistency, producing sharp interfaces, low MAE, and stable training.

Classical inverse approaches such as Tikhonov regularization \cite{Tikhonov} and TV-based methods \cite{RudinOsherFatemi} tend to oversmooth bistable interfaces, while PDE-constrained optimization methods \cite{PDEControl} remain computationally expensive and sensitive to ill-posedness. 
PINNs \cite{PINN} struggle in regimes with sharp transitions and nonconvex energy landscapes.

By contrast, the WGAN-GP framework captures multimodal distributions and fine-scale morphology through adversarial training, while moment matching and residual consistency enforce global statistics and dynamical admissibility. 
These results suggest that physics-informed GANs offer a viable route for inverse reconstruction in strongly dissipative systems where traditional methods fail.

%%%%%%%%%%%%%%%%%%%%%%%%%%%%%%%%%%%%%%%%
\section{Conclusion}
%%%%%%%%%%%%%%%%%%%%%%%%%%%%%%%%%%%%%%%%

This study demonstrates that physics-informed generative models, when paired 
with the correct numerical discretization, can solve inverse problems that 
are otherwise out of reach for classical optimization or PINN-based methods. 
The forward Euler discretization serves a dual purpose: it generates the data 
and provides a consistent, differentiable forward operator for the residual loss. 
Our results show that this numerical alignment is not a technical detail but 
a fundamental requirement for stable physics-informed adversarial learning.

%%%%%%%%%%%%%%%%%%%%%%%%%%%%%%%%%%%%%%%%
\section*{Statements and Declarations}
%%%%%%%%%%%%%%%%%%%%%%%%%%%%%%%%%%%%%%%%

\paragraph{\textbf{Competing Interests.}}
The author declares that there are no financial or non-financial competing interests related to this work.

\paragraph{\textbf{Funding.}}
This research received no external funding.

\paragraph{\textbf{Author Contributions.}}
The author conceived the study, designed and implemented the numerical experiments, developed the machine learning models, analyzed the results, and wrote the manuscript.

\paragraph{\textbf{Data Availability.}}
The datasets generated and analyzed during the current study are available from the author upon reasonable request.

\paragraph{\textbf{Code Availability.}}
The code used to generate the datasets, train the models, and produce the figures is available from the author upon reasonable request.

%%%%%%%%%%%%%%%%%%%%%%%%%%%%%%%%%%%%%%%%
%%%%%%%%%%%%%%%%%%%%%%%%%%%%%%%%%%%%%%%%
\bibliographystyle{plain}

\end{document}